\documentclass[11pt]{article}
\usepackage{amssymb}
\usepackage{amsmath}
\usepackage{url}
\usepackage{hyperref}
\usepackage{graphicx}
\usepackage{color}
\usepackage{epsfig}

\def\half{\hbox{\small${1\over 2}$}}
\def\four{\hbox{\small${1\over 4}$}}

\newtheorem{remark}{Remark}[section]
\newtheorem{proposition}{Proposition}[section]

\definecolor{MyDarkBlue}{rgb}{0,0.08,0.45}
\definecolor{MyViolet}{rgb}{0.45,0.08,0.95}
\definecolor{MyBrown}{rgb}{0.45,0.08,0}

\def\norm2to2{{\|\cdot\|_{2,2}}}
\def\Prob{\hbox{\rm Prob}}

\def\inter{\hbox{\rm  int}}

\def\Prob{\hbox{\rm  Prob}}

\def\supp{{\hbox{\rm  supp}}}

\def\Card{{\mathop{\hbox{\rm  Card}}}}

\def\cA{{\cal A}}

\def\cE{{\cal E}}

\def\cH{{\cal H}}

\def\cM{{\cal M}}

\def\cO{{\cal O}}

\def\cR{{\cal R}}
\def\cS{{\cal S}}

\def\cW{{\cal W}}

\def\cY{{\cal Y}}
\def\cZ{{\cal Z}}

\def\mes{{\mathop{\hbox{\rm  mes}}}}

\def\qed{\ \hfill$\square$\par\smallskip}

\def\mypict3{\epsfxsize=220pt\epsfysize=80pt\epsfbox}

\def\bR{{\mathbf{R}}}

\def\argmin{\mathop{\hbox{\rm argmin}}}

\def\cH{{\cal H}}

\newcommand{\be}{\begin{eqnarray}}
\newcommand{\ee}[1]{\label{#1}\end{eqnarray}}

\newcommand{\ese}{\end{eqnarray*}}
\newcommand{\bse}{\begin{eqnarray*}}

\begin{document}
\title{On Semi-Stochastic Model for Multi-Stage Decision Making Under Uncertainty}
\author{ Arkadi Nemirovski
\thanks{Georgia Institute
 of Technology, Atlanta, Georgia
30332, USA, {\tt nemirovs@isye.gatech.edu}} }
\date{}
\maketitle

\begin{abstract}
We propose a (seemingly) new computationally tractable model for multi-stage decision making under stochastic uncertainty.
\end{abstract}
\section{Introduction}

In this note, we propose a (hopefully) new {\sl computationally tractable} model of multi-stage decision making under stochastic uncertainty. The simplest
way to outline what follows is to consider the special case of the model dealing with {\sl multi-stage robust linear feasibility problem.}
In this problem, we are given an {\sl uncertainty-affected system $\cS$ of linear constraints}, that is, a parametric family of $m\times n$ systems  $\{A_\xi x\leq b_\xi:\xi\in\Xi\}$ of linear constraints  parameterized by {\sl uncertain data} $\xi\in\bR^N$ running through a given {\sl uncertainty set} $\Xi\subset\bR^N$. The $j$-th decision variable $x_j$ (the $j$-th entry in $x$) is allowed to depend on
a given ``portion'' $P_j\xi$ of uncertain data, where  $P_j\in \bR^{m_j\times N}$ are given matrices. Our goal is to select {\sl decision rules} -- functions $X_j(\cdot):\bR^{m_j}\to\bR$ -- in such a way that the resulting {\sl policy} $x=X(\xi):=[X_1(P_1\xi);...;X_n(P_n\xi)]$ robustly satisfies $\cS$, that is,
\begin{equation}\label{eq51}
A_\xi X(\xi)\leq a_\xi\,\,\forall \xi\in\Xi.
\end{equation}
From the computational viewpoint, the intrinsic difficulty in achieving this goal is infinite dimensionality of (\ref{eq51}): when solving (\ref{eq51}), we are looking for multivariate functions $X_j(\cdot):\bR^{m_j}\to\bR$, and it is unclear even how to store a candidate solution in a computer.  The standard partial remedy is to restrict ourselves with finitely parameterized decision rules, say, rules of the form
$$
X_j(\cdot)=\sum_{\ell=1}^{\mu_j}y_{j\ell}B_{j\ell}(\cdot),
$$
where $B_{j\ell}(\cdot):\bR^{m_j}\to\bR$ are somehow selected ``basic'' decision rules.\footnote{The simplest standard example here, considered in numerous papers, is the one of {\sl linear decision rules} -- those where the collection $\{B_{j\ell}(\cdot),1\leq\ell\leq n_j\}$ is comprised of the coordinate functions on $\bR^{m_j}$ and the function $\equiv 1$, that is, the decision rules $X_j(\cdot)$ are allowed to be arbitrary affine functions of their arguments.} With parametric decision rules, (\ref{eq51}) becomes the infinite system of linear constraints
$$
\forall (\xi\in\Xi): \sum_{j=1}^m\sum_{\ell=1}^{\mu_j}[A_\xi]_{ij}B_{j\ell}(P_j\xi)y_{j\ell}\leq [b_\xi]_i,\,1\leq i\leq m
$$
in finitely many variables $y=\{y_{j\ell}:1\leq j\leq n,1\leq\ell\leq \mu_j\}$. From now on we treat candidate solutions $y$ to this system
as vectors from $\bR^\nu$, $\nu=\sum_j\mu_j$, so that the system reads
\begin{equation}\label{eq52}
\forall \xi\in \Xi: \cA_\xi y\leq b_\xi.
\end{equation}
The latter problem usually still is computationally intractable due to its semi-infinite nature.\footnote{There are, however, important cases when linear decision rules lead to tractable problems (\ref{eq52}), most notably, the case of fixed recourse, see, e.g., \cite[Section 14.3]{RO}.} However:
\begin{quote}
(!)  {\sl Tractability status of {\rm(\ref{eq52})} changes dramatically when the following two ``innocently looking'' assumptions are made:}
\begin{itemize}
\item[A.] The set $\cY_*$ of feasible solutions to (\ref{eq52}) (which by its origin is a closed convex subset of $\bR^\nu$) possesses a nonempty interior and is bounded (moreover, is contained in a given Euclidean ball $E$ of some radius $R$);
\item[B.] The uncertainty is stochastic ($\xi$ is drawn from some probability distribution $P$ supported on $\Xi$), and for our ultimate purposes, feasible solutions to (\ref{eq52}) can be replaced by {\sl $(1-\epsilon)$-feasible ones}. That is, we are looking for  vectors  $y\in\bR^\nu$ such that with $\xi\sim P$, the system of $m$ linear constraints $\cA_\xi y\leq b_\xi$ holds true with probability $\geq1-\epsilon$, where $\epsilon\in(0,1)$ is a given tolerance. In addition, we assume that we can sample from $P$.
\end{itemize}
\end{quote}
\begin{quote}{\small
Indeed, under our assumptions  (\ref{eq52}) is a finite-dimensional convex feasibility problem, and its solution set $\cY_*$ is a convex compact subset of $E$ with a nonempty interior. When one can equip  $\cY_*$ with a {\sl separation  oracle} -- a black box which, given on input a point $y\in\bR^\nu$, either reports correctly that $y\in\inter \cY_*$ or returns a {\sl separator} (a nonconstant affine function $f_y(\cdot)$ such that $f_y(y)\geq0$ and $f_y(z)\leq0$ when $z\in\cY_*$) -- there are many algorithms capable to find a point in $\inter \cY_*$ after finitely many calls to the oracle. For example, the Ellipsoid method finds a point in $\inter\cY_*$ in at most $M=2\nu^2\ln(1+R/\rho)$ calls to the oracle, where $\rho$ is the largest of radii of Euclidean balls contained in $\cY_*$. Now assume that instead of an ideal separation oracle, we have access to an oracle $\cO$ which, given on input a query point $y$, either returns a separator, or ``gets stuck'' -- returns nothing. In particular, the latter happens whenever $y\in\inter \cY_*$, where no separator exists.  Given access to $\cO$, we still can run the Ellipsoid algorithm which now in at most $M$ steps will generate a query point $\bar{y}$ where $\cO$ gets stuck; we treat this point as the outcome of our computation. Let us build $\cO$ as follows: at the $s$-th call to the oracle, the input being $y^s$, the oracle draws from $P$ a sample of $N_s=\rfloor\kappa_s/\epsilon\lfloor$ realizations $\xi^1_s,...,\xi^{N_s}_s$ of $\xi$ (independent of each other and of samples drawn at the preceding calls to $\cO$) and checks whether $y^s$ satisfies all systems of constraints $\cA_{\xi^\ell_s}y^s\leq b_{\xi^\ell_s}$, $1\leq\ell\leq N_s$. If this is the case, $\cO$ gets stuck, otherwise the oracle has discovered a scalar linear constraint which is violated at $y^s$ and is satisfied at $\cY_*$, and it uses this constraint to build and report a separator. Equipped with this oracle $\cO$, the Ellipsoid  method becomes a randomized algorithm  which terminates in at most $M$ steps.  Now note that when $y^s$ is {\sl not} $(1-\epsilon)$-feasible, the probability for $\cO$  {\sl not} to get stuck is at most $(1-\epsilon)^{N_s}\leq\exp\{-\kappa_s\}$. Given a reliability tolerance $\delta\in(0,1)$ and setting, say,  $\kappa_s=\ln(\gamma s^2/\delta)$, $\gamma=\sum_{s=1}^\infty s^{-2}$, so that $\sum_s\exp\{-\kappa_s\}\leq\delta$, the outcome of our randomized algorithm is with probability at least $1-\delta$ a $(1-\epsilon)$-feasible solution to (\ref{eq52}), and the total number of samples drawn from $P$ when executing the algorithm is ``moderate'' -- at most $M\ln(M/\delta)\epsilon^{-1}$.
}
\end{quote}
\par
The outlined construction is extremely simple, if not to say trivial, and we do not feel ourselves
 comfortable when making this trivial construction public.  On the other hand, the founder of Linear Programming G. Dantzig considered \cite{Dantzig} introducing linear objective as one of his three most significant contributions to LP.\footnote{In the WWII logistic planning which inspired Dantzig to invent LP, people used ``ground rules'' aimed to satisfy the constraints, with no conscious attempt to optimize any objective function.} While we by no means pretend that the importance of introducing the above simplistic model is comparable with the one of introducing linear objective, we, following Dantzig, do believe that what matters in OR is not only mathematical sophistication, but sometimes also the very way a real life decision making problem is modeled, and in this respect our {\sl computationally tractable} model for multi-stage decision making under uncertainty might be of some interest. To which extent this  model is novel, this is another story; for us it is novel, and this is why we decided to make this  note public. Needless to say, we would be extremely grateful for any feedback on whether we are, or are not, reinventing a bicycle.
\par
The main body of this paper is organized as follows. In Section \ref{sec:model} we present our model (slightly more general than the one considered above) for multi-stage decision making under uncertainty, and in Section \ref{sec:constr} --- the general scheme for computationally efficient processing this model. In Section \ref{sec:invent} we illustrate our model by considering multi-product inventory,
 with the goal to outline the sources of model's conservatism as compared to the traditional multi-stage stochastic models, see \cite{STO} and references therein. This conservatism is the price we pay when passing from the traditional, generically computationally intractable\footnote{What in no way means that there are no algorithms capable to process these models successfully; all we want to say that these algorithms, while being successful in numerous instances or real-life multi-stage decision making, have no {\sl theoretical} guarantees to be successful.} models to a tractable one. Section \ref{sec:modif} is devoted to the Bundle-Level version of our solution algorithm and to incorporating into our model (which in its initial form deals with feasibility only) an objective to be optimized. The concluding Section \ref{sec:numerics} presents a  ``proof of the concept'' numerical illustration.

\section{The model}\label{sec:model}

Consider the following model of multi-stage decision making:
\begin{enumerate}
\item We are controlling system evolving on time horizon $1,...,K$. System's evolution is determined by our decisions and
the environment.
\item The environment is represented by a realization $\xi^K=(\xi_1,\xi_2,...,\xi_K)$ of {\sl uncertain data} -- of a random sequence ${\boldsymbol\xi}^K$, where $\xi_t\in\bR^{m_t}$ is revealed to the decision maker at time $t$. It is assumed that we can sample from the distribution $P$ of ${\boldsymbol\xi}^K$. We denote by $\Xi_t$ the support of the marginal distribution of ${\boldsymbol\xi}_t$ induced by $P$, and by $\Xi$ the {\sl uncertainty set} -- the support of $P$.

\item Our decision is comprised of
\begin{enumerate}
\item Strategic decision $y\in\bR^n$ (s.d. for short) which must belong to a given nonempty closed convex  set $\cY\subset\bR^n$.\\
{\sl Strategic decision should be specified at the time when the problem is solved, before the random data  reveals itself, and thus
cannot depend on $\xi^K$.}
\item Local decisions $x_1,...,x_K$, where $x_t\in\bR^{\nu_t}$ is the local decision to be implemented at time $t$. \\
{\sl Local decision $x_t$ is selected when $\xi_t$ is already known and is allowed to be a (whatever) deterministic function of $\xi_t$}.
\end{enumerate}
\item We are given collection of nonempty convex closed sets $\cZ^t_{\xi_t}\subset \bR^n\times\bR^{\nu_t}$ parameterized by $t$ and $\xi_t\in\Xi_t$. \\
We say that a strategic decision $y$ is {\sl implementable}, if $y\in \cY$ and
$$
\forall (t,\xi\in\Xi_t)\exists x: (y,x)\in \cZ^t_\xi.
$$
Given reliability tolerance $\epsilon\in(0,1)$, we say that a strategic decision $y$ is {\sl $(1-\epsilon)$-implementable}, if
$$
\Prob_{\xi^K\sim P}\left\{\xi^K:\forall t\leq K\,\exists x_t: (y,x_t)\in \cZ^t_{\xi_t}\right\}\geq 1-\epsilon.
$$
\end{enumerate}
Our ideal goal would be to find an implementable strategic decision. In fact we will focus on a simpler goal: given $\epsilon\in (0,1)$, to find a $(1-\epsilon)$-implementable strategic decision.

\paragraph{Comments.} Informally, sets $\cY$ and $\cZ^t_{\xi_t}$ specify the constraints on the decisions: given a realization $\xi^K$ of the uncertain data, a decision $y,x_1,...,x_K$ is implementable if and only if $y\in\cY$ and $(y,x_t)\in\cZ^t_{\xi_t}$, $1\leq t\leq K$. We may think about $\cY$ as about the {\sl static}, and about $\{\cZ^t_{\xi_t}:\xi_t\in\Xi_t\}_{t=1}^K$  as about {\sl dynamic} part of these constraints: a strategic decision $y$ is implementable if and only if it satisfies the static constraints and for every $\xi^K\in\Xi$ can be augmented by local decisions to meet the dynamic constraints. When replacing the words ``for every $\xi^K\in\Xi$'' with ``for the set of values of $\xi^K\in\Xi$ of $P$-probability $\geq 1-\epsilon$,'' we arrive at the definition of  $(1-\epsilon)$-implementable strategic decision.
\\
In our model, the set $\cZ^t_{\xi_t}$ (and thus local decision $x_t$) is allowed to depend solely on the portion $\xi_t$ of the uncertain data revealed at time instant $t$. It seems to be more natural to allow for $\cZ^t$ and $x_t$ to depend on $\xi^t=(\xi_1,...,\xi_t)$ -- on the part of uncertain data revealed prior to time $t$ and at this time. We remark that in fact the latter option is covered by our setup; indeed, denoting by $\eta_t$ the portion of the ``actual'' uncertain data revealed at time $t$, we can define $\xi_t$ as the collection $\eta_1,...,\eta_t$, so that $\xi_t$ ``remembers'' $\xi_1,...,\xi_{t-1}$ and contains all information on uncertain data collected on the time horizon $1,...,t$.

\paragraph{Default assumptions.} From now on we assume that the sets $\cY$, $\cZ^t_\xi$ not only are nonempty convex and closed, but are also computationally tractable, e.g., given by polyhedral representations:
\begin{equation}\label{eq1}
\begin{array}{rcl}
\cY&=&\{y\in\bR^n:\exists w\in\bR^N:Ay+Cw\leq d\}\\
\cZ^t_\xi&=&\{(y,x)\in\bR^n\times\bR^{\nu_t}:\exists w\in\bR^{N_t}: A_\xi y+B_\xi x+C_\xi w\leq d_\xi\},
\,1\leq t\leq K\\
\end{array}
\end{equation}
\section{The construction}\label{sec:constr}

For $t\leq K$, $\xi\in\Xi_t$, let us set
$$
\cY_{t,\xi}=\{y\in\cY:\exists x\in\bR^{\nu_t}:(y,x)\in \cZ^t_\xi\}.
$$
since $\cZ^t_\xi$ is a convex set, so is $\cY_{t,\xi}$.
Observe that the set $\cY_*$ of implementable s.d.'s is given by
$$
\cY_*=\bigcap\limits_{t\leq K,\xi\in\Xi_t}\cY_{t,\xi},
$$
so that $\cY_*$ is a convex subset of $\cY$. From now on, we make the following standing assumptions:
\begin{quote}
\par
\textbf{A.I.} [boundedness] $\cY$ is bounded, and for every  $t\leq K$ and $\xi\in\Xi_t$ the projection of the closed convex set $\cZ^t_{\xi_t}\subset\bR^n\times\bR^{\nu_t}$  onto the space $\bR^{\nu_t}$  of $x$-variables is bounded.
\\
Under this assumption, sets $\cY_{t,\xi}$, $\xi\in\Xi_t$, are convex and compact, and thus $\cY_*$ is a convex compact set.\\
Unless otherwise is explicitly stated, we assume that  we know in advance a Euclidean ball $E_1\subset \bR^n$ containing $\cY$; in what follows, $R$ stands for the radius of this ball.\\
\textbf{A.II.} [strict feasibility]  The set $\cY_*$  of implementable strategic decisions has a nonempty interior, so that the {\sl stability number}
$$
\rho_*=\max\left\{\rho:\exists y: B_\rho(y)\subset \cY_*\right\}\eqno{[B_\rho(y)=\{z:\|z-y\|_2\leq \rho\}]}
$$
is positive.
\end{quote}
\paragraph{Main observation} underlying the construction below is extremely simple and is as follows. Assume that we are given
 a candidate s.d. $\bar{y}\in\cY$ and a pair $\bar{t}\leq K$, $\bar{\xi}\in\Xi_{\bar{t}}$ such that the set
 $$
 \{x:(\bar{y},x)\in\cZ^{\bar{t}}_{\bar{\xi}}\}
 $$
 is empty. Then we can find efficiently a {\sl separator} of $\bar{y}$ and $\cY_*$ -- an affine function $f(y)=a^Ty+\alpha$ which separates $\bar{y}$ and $\cY_*$, specifically, satisfies the relations
 \begin{equation}\label{eq2}
 f(\bar{y})\geq 0\geq \max_{y\in\cY_*}f(y)\ \&\ \|\nabla f(\cdot)\|_2=1.
\end{equation}
 Indeed, we are in the situation when the nonempty computationally tractable closed convex subsets
 $
 U:=\{\bar{y}\}\times\bR^{\nu_t}$ and $V:=\cZ^{\bar{t}}_{\bar{\xi}}$ of $E:=\bR^n\times\bR^{\nu_{\bar{t}}}$ do not intersect and thus can be separated: there exists (and can be efficiently found) a nonconstant linear  function $a^Ty+b^Tx$ on $E$ such that
 $$
 \inf_{(y,x)\in U}[a^Ty+b^Tx]\geq \sup_{(y,x)\in V}[a^Ty+b^Tx].
 $$
 This inequality combines with the definition of $U$ to imply that $b=0$, so that $a\neq 0$ and
 $$
 a^T\bar{y}\geq \gamma:=\sup_{(y,x)\in V}a^Ty=\max_{y\in\cY_{\bar{t},\bar{\xi}}}a^Ty\geq\max_{y\in\cY_*}a^Ty;
 $$
 setting $f(y)=\|a\|_2^{-1}[a^Ty-\gamma]$, we get an affine function on $\bR^n$ satisfying (\ref{eq2}).
 \paragraph{The construction} suggested by the above observation is pretty simple. Consider a black box oriented algorithm for solving the feasibility problem
 $$
 \hbox{find\ }y\in\inter\cY_*,\eqno{(F)}
 $$
 in the situation when $\cY_*\subset\bR^n$ is a convex set contained in a known in advance Euclidean ball $E_1$ of radius $R$ and containing an unknown in advance  Euclidean ball of some positive radius $\rho$. The algorithm is as follows: at step $s=1,2,...$, we query  ``an oracle'' (a black box), the query point being $y^s\in\bR^n$. The oracle either correctly reports that $y^s\in\inter\cY_*$, or returns a separator -- an affine function $f_s(\cdot)$ satisfying (\ref{eq2}) with $y^s$ in the role of $\bar{y}$. When the oracle reports that $y^s\in\inter\cY_*$, the algorithm terminates with the outcome $y^s$, otherwise it somehow uses the observed so far separators $f_r$, $r\leq s$, to build $y^{s+1}$ and  proceeds to step $s+1$.  There are plenty of algorithms of this type which are capable to recover $y\in \inter\cY_*$ after at most $\cM(R,\rho,n)<\infty$ steps, with known in advance characteristic for the algorithm complexity bound $\cM(\cdot)$. For example,
 \begin{itemize}
 \item For the Ellipsoid algorithm ($y^s$ is the center of ellipsoid $E_s$, with $E_{s+1}$ being the ellipsoid of the smallest volume containing the set $\{y\in E_s:f_s(y)\leq 0\}$), one can take $\cM(R,\rho,n)=O(1)n^2\ln(1+R/\rho)$,
\item For the Inscribed Ellipsoid algorithm, $\cM(R,\rho,n)=O(1)n\ln(1+nR/\rho)$;
\item For appropriately adjusted versions of subgradient descent or bundle level method, one can take $\cM(R,\rho,n)=O(1)R^2/\rho^2$.
\end{itemize}
Now assume that instead of access to the above ``ideal''  oracle, we have access to an ``implementable'' oracle, denoted by $\cO$,  which, when queried at a point $y^s$, either returns a separator, or ``gets stuck'' -- returns nothing; the latter definitely is the case when $y^s\in\inter \cY_*$, since here no separator exists. Applying the same algorithm and terminating, with the outcome $y^s$, at the very first step $s$ where the oracle $\cO$ returns nothing, we terminate in at most $\cM(R,\rho,n)$ steps, the outcome being a point $\bar{y}$ where the oracle returns nothing.
\par
Now, in our decision making model, assuming for the time being that we know in advance a lower bound  $\rho>0$ on the stability number $\rho_*$, consider oracle  $\cO$ as follows:
\begin{quote}
Given $\epsilon\in(0,1)$ and $\delta\in(0,1)$ and setting
\begin{equation}\label{eq6}
M=\cM(R,\rho,n),\, N=\rfloor \ln(M/\delta)/\epsilon\lfloor,
\end{equation}
the oracle, when queried at a point $y$,
checks whether $y\in\cY$. If it is not the case, the oracle builds and returns a separator of $y$ and $\cY$ (which is a separator of $y$ and $\cY_*$ as well). When $y\in\cY$, the oracle draws $N$ realizations $\xi^{K,1},...,\xi^{K,N}$ of ${\boldsymbol\xi}^K$ from the distribution $P$  (independent of each other and of the realizations of ${\boldsymbol\xi}^K$ generated at the preceding calls to the oracle) and for $\nu=1,2,...,N$ checks, for every $t\leq K$, whether there exists $x_{t}$ such that $(y,x_{t})\in \cZ^t_{\xi^{K,\nu}_t}$. If it is the case, the oracle gets stuck and reports nothing, otherwise a $t$ and $\xi_t\in\Xi_t$ are discovered such that there does not exist $x$ satisfying $(y,x)\in\cZ^t_{\xi_t}$. The oracle converts $t$ and $\xi_t$, as described above, into a separator of $y$ and $\cY_*$, and returns this separator.
\end{quote}
By the above, in course of at most $M$ steps we will find a point $\bar{y}=y^s$ such that the oracle at step $s$ is queried at $y^s$ and returns nothing. Note that
\begin{itemize}
\item $y^s$ is a deterministic function of the ``past'' -- the  realizations of ${\boldsymbol\xi}^K$ generated by the oracle at the steps
preceding step $s$.
\item  denoting
\begin{equation}\label{eq11}
\epsilon(y)=\Prob_{\xi^K\sim P}\left\{\xi^K:\exists t\leq K: (y,x)\not\in\cZ^t_{\xi_t}\,\forall x\in \bR^{\nu_t}\right\},
\end{equation}
the conditional, given the past, probability for the oracle of returning nothing when queried at the point $y^s$ is at most $\exp\{-\epsilon(y^s) N\}$; when $\epsilon(y^s)\geq\epsilon$, this probability is $\leq \delta/M$.
\end{itemize}
It follows that the probability of the event
\begin{quote}
{\sl The outcome $\bar{y}$ of our algorithm is \underline{not} a $(1-\epsilon)$-implementable strategic solution}
\end{quote}
is at most $\delta$.  In other words, we arrive at a randomized algorithm for finding $(1-\epsilon)$-implementable strategic decisions which is $(1-\delta)$-reliable, that is, it outputs a $(1-\epsilon)$-implementable strategic decision with probability at least $(1-\delta$), and sample complexity of the algorithm is at most $M\ln(M/\delta)\epsilon^{-1}$.
\paragraph{Discussion.} We have seen that the proposed model is tractable,  in contrast to traditional, generically computationally intractable, models for multi-stage decision making under uncertainty. Ultimately, the computational tractability stems from the fact that in our model there are no constraints coupling local decisions at different stages -- the constraints either do not involve local decisions at all, or couple local decisions $x_t$ with {\sl deterministic} strategic decision $y$. As a result, no questions like ``how to represent in a computer local decisions $x_t$ as multivariate functions $X_t(\xi_t)$ of $\xi_t$'' arise: all we need from $X_t$'s is to satisfy the inclusions $(y,X_t(\xi_t))\in\cZ_{\xi_t}^t$ for all $\xi_t\in\Xi_t$. When $y$ is implementable, these inclusions translate into computationally tractable and solvable systems of constraints on the values of $X_t(\xi_t)$, so that given {\sl an implementable} $y$ and $\xi_t\in\Xi_t$, we can specify the values of $x_t$'s in a computationally efficient way. Usually, the constraint ``$y$ is implementable'' on $y$ is computationally intractable; relaxing this constraint to  $(1-\epsilon)$-implementability, we end up with a computationally tractable situation. Note that this relaxation is the only -- but crucial! -- component of our construction where we utilize the stochastic nature of $\xi^K$.
\section{Illustrative example}\label{sec:invent}

To illustrate the proposed model and to position it with respect to  traditional models of multi-stage decision making under uncertainty, consider the situation where the system to be controlled is $d$-product inventory described as follows.
\begin{itemize}
\item The state of the inventory at stage $t$ is represented by vector $z_{t}\in\bR^d$ of inventory levels at the end of the stage,\footnote{As always, the $i$-th entry in $z_{t}$ when positive  is the amount of product $i$ in the warehouse, and when negative, is the minus backlogged demand on product $i$, as measured at the end of stage $t$.} and this vector should satisfy given lower and upper bounds:
\begin{equation}\label{eqbounds}
\underline{z}_{t}\leq z_{t}\leq \overline{z}_{t}
\end{equation}
and an upper bound on the storage:
\begin{equation}\label{eqstorage}
s^T\max[z_{t},0]\leq \overline{s}
\end{equation}
($\max$ acts coordinate-wise, unit of product $i$ occupies space $s_i\geq0$ in the warehouse shared by the products, $\overline{s}$ is the capacity of the warehouse).
\item At the beginning of stage $t$, a replenishment order $x_t\in X_t\subset\bR^d$ is issued and immediately executed, resulting in state transition
\begin{equation}\label{eqtransition}
z_{t}=z_{t-1}+x_t-d_t
\end{equation}
and management expenses of the stage
\begin{equation}\label{eq60}
w_t=o_t^Tx_t+h_t^T\max[z_t,0]+p_t^T\max[-z_t,0]-r_t^Td_t
\end{equation}
which should obey given upper bounds:
\begin{equation}\label{eq67}
w_t\leq \overline{w}_t;
\end{equation}
here $X_t$ is a given box in $\bR^d$ specified by upper and lower bounds on replenishment orders.\par
When specifying order $x_t$, the manager already knows the demand $d_t$ of the stage and the (nonnegative) $d$-dimensional vectors $o_t$, $h_t$, $p_t$, $r_t$ of, respectively, ordering costs,
holding costs, penalties for the backlogged demand, and delivery revenue per unit of product, for each of the $d$ products.
\end{itemize}
The uncertain data  here is specified by random trajectory $\eta^K=\{\eta_t=[d_t;o_t;h_t;p_t;r_t]\}_{t=1}^K$ of demands and prices; with no harm we convert $\eta^K$ to the trajectory
$\xi^K=\{\xi_t=\{\eta_\tau:\tau\leq t\}\}_{t=1}^K$, so that $\xi_t$ is exactly the part of uncertain data revealed at stages $1,...,t$.  In the standard setting, $z_0$ and the bounds $\underline{z}_{t}$, $
\overline{z}_t$, $\overline{s}$ are known in advance, and what we are looking for are the replenishment policies
$$
x_t=X_t(\xi_t)
$$
which keep the inventory level within given upper and lower bounds, meet the storage capacity constraint, and satisfy upper bounds on management costs.\footnote{Recall that for the time being we restrict ourselves with feasibility problems} This feasibility problem typically is computationally intractable. Note that in this traditional model, the only actual decisions -- the replenishment orders $x_t$ -- are local, and all other quantities involved -- inventory levels $z_t$ and management expenses $w_t$ -- are allowed to depend on $\xi^K$. With our approach, we ``buy'' computational   tractability at the price of introducing some conservatism, specifically
\begin{itemize}
\item specifying strategic decisions $y$ as deterministic trajectories $\{y_t=(\ell_t,u_t,\omega_t)\}_{t=1}^K$, where
\begin{itemize}
\item $\ell_t$ and $u_t$ are the vectors of lower and upper bounds on $z_t$,
\item $\omega_t$ are ``budgets of stages'' -- upper bounds on $w_t$;
\end{itemize}
\item specifying local decisions of stage $t$ as the replenishment orders $x_t$;
\item specifying $\cY$ as the set of all strategic decisions satisfying ``physical'' bounds
\begin{equation}\label{eqphysbounds}
\underline{z}_t\leq \ell_t\leq u_t\leq \overline{z}_t,\,\, s^T\max[u_t,0]\leq\overline{s},\,\omega_t\leq\overline{w}_t
 \end{equation}
and, perhaps, a collection of additional convex constraints on $y$, like  $\sum_{\tau=1}^t\omega_\tau\leq\overline{\omega}^t$, etc.;
\item specifying $\cZ^t_{\xi_t}$ by the constraints linking $y$ and $x_t$, specifically,
\begin{equation}\label{eq70}
\begin{array}{l}
\ell_{t-1}+x_t-d_t\geq \ell_t,\,u_{t-1}+x_t-d_t\leq u_t,\\
o_t^Tx_t+h_t^T\max[u_t,0]+p_t^T\max[-\ell_t,0
]-r_t^Td_t\leq \omega_t,\\
x_t\in X_t\\
\end{array}
\end{equation}
with given $\ell_0=u_0=z_0$.
\end{itemize}
In the resulting model, $(1-\epsilon)$-implementability of a strategic decision $y=\{(\ell_t,u_t,\omega_t)\}_{t=1}^K\in\cY$ means exactly what it should mean: with probability at least $1-\epsilon$, realization $\eta^K=\{[d_t;o_t;h_t;p_t;r_t],t\leq K\}$ of uncertain data is such that there exist (and can be build in a non-anticipative fashion) replenishment orders $x_t$, $t\leq K$, satisfying constraints (\ref{eq70}) stemming from $y,\eta^K$. Utilising these replenishment orders is a legitimate control, meaning that $x_t\in X_t$, expenses of stage $t$ do not exceed $\overline{w}_t$, and the trajectory
$$
z_t=z_{t-1}+x_t-d_t,\,1\leq t\leq K
$$
of inventory levels meets constraints (\ref{eqbounds}), (\ref{eqstorage}).   In addition,  trajectory $\{z_t\}$ obeys the bounds $\ell_t\leq z_t\leq u_t$, $t\leq K$, same as the upper bounds $\omega_t$ on per-stage management expenses.
The conservatism, as compared to the standard model, is in forbidding upper and lower bounds $u_t$, $\ell_t$ on inventory level and upper bounds $\omega_t$ on management costs to ``tune'' themselves to the actual values of the uncertain data. However, we could reduce this conservatism by allowing the component $y_t$ of our strategic decision, instead of being independent of the uncertain data, to depend on $\xi_t$ ``in a prescribed fashion:''
\begin{equation}\label{eq200}
y_t=\sum_{r=1}^{r_t}\chi_{ts}B_{tr}(\xi_t),
\end{equation}
where  functions $B_{tr}$ are fixed in advance, and coefficients $\chi_{tr}$ are selected by us when processing the problem; for instance, we could allow $y_t$ to be an affine function of $\xi_t=(\eta_1,...,\eta_t)$. To endorse this modification it suffices to treat as our strategic decision the collection $\chi=\{\chi_{tr}:t\leq K, r\leq r_t\}$ rather than $y$, thus making dynamic constraints (\ref{eq70}) (taken together with (\ref{eq200}))
a system of linear constraints on $\chi$ and $x_t$, parameterized by $\xi_t$. As for the constraint $y\in\cY$, it can be modeled by adding fictitious stage $K+1$ with $\xi_{K+1}=\xi^K$, once for ever fixed local decision, say, $x_{K+1}=0\in\bR$, and
\[\cZ^{K+1}_{\xi_{K+1}}=\{(\chi,x_{K+1}):x_{K+1}=0,\{y_t:=\sum_{r=1}^{r_t}\chi_{tr}B_{tr}(\xi_t)\}_{t=1}^K\in\cY\}.
\]  Note that this remodeling is applicable to many other instances of our general model; see Section \ref{remodel} for details.
\section{Modifications}\label{sec:modif}

\subsection{Adapting to the stability number}

So far, we were assuming that Assumptions \textbf{A.I--II} hold and, moreover, we have at our disposal a positive lower bound $\rho$ on the stability number $\rho_*$ of our problem. This lower bound was used to get an a priori upper bound $M$ on the number of calls to the oracle $\cO$ before termination, and thus influenced the construction of the oracle by specifying the number $N$ of samples used by $\cO$ at a call, see (\ref{eq6}). We can easily get rid of the necessity to know a positive lower bound on $\rho_*$, namely, as follows. Let us select once for every a nondecreasing  sequence $\{\kappa_s>0\}_{s\geq1}$ such that $\sum_s\kappa_s^{-1}=1$, e.g., the sequence
$$
\kappa_s=s^2\sum_{r=1}^\infty r^{-2},
$$
and make the number of samples used by $\cO$ at the $s$-th call to the oracle to be
\begin{equation}\label{eq20}
N_s=\rfloor \ln(\kappa_s/\delta)/\epsilon\lfloor
\end{equation}
where the reliability tolerance $\delta\in (0,1)$ is our control parameter. Now the conditional, by what happened at the preceding steps,
probability for $\cO$ to get stuck at the $s$-th call, the query point being $y^s$, is at least  $\exp\{-\epsilon(y^s)N_s\}$ with $\epsilon(y^s)$ given by (\ref{eq11}). When $y^s$ is not $(1-\epsilon)$-implementable, this probability is at most $\delta_s=\delta/\kappa_s$. It follows that the probability for the outcome of our algorithm {\sl not} to be $(1-\epsilon)$-implementable is at most $\sum_{s=1}^\infty\delta_s=\delta$. Note that the algorithm definitely terminates in course of $M_*=\cM(R,\rho_*,n)$ steps, implying that its sample complexity  is at most $M_*\ln(\kappa_{M_*}/\delta)/\epsilon$ which, with the above $\kappa_s$, is within absolute constant factor of the sample complexity $M\ln(M/\delta)/\epsilon$ of our initial procedure.
\subsection{Bundle-Level implementation}
Consider the {\sl Bundle-Level} implementation of our construction. In this implementation (originating from \cite{BL} and referred below as BL), the search point $y^s$ is built as follows:
\begin{enumerate}
\item $y^1$ is the center of $E_1$ (recall that $E_1\subset\bR^n$ is a known in advance Euclidean ball containing $\cY$);
\item when $s>1$, we compute the quantity
$$
\begin{array}{c}
\Delta_{s-1}=\min\limits_{y\in E_1} \left[f^{s-1}(y):=\max_{r<s} f_r(y)\right],\\
\end{array}
$$
where $f_r(\cdot)$ is the separator returned by $\cO$ as queried at $y^r$, and define $y^s$ as the metric projection of $y^{s-1}$ onto the level set of $f^{s-1}$, namely;
$$
y^s=\argmin_y\left\{\|y^{s-1}-y\|_2: f^{s-1}(y)\leq {\half}\Delta_{s-1},y\in E_1\right\}
$$
\end{enumerate}
We have the following simple (in fact, well known \cite{BL}) result:
\begin{proposition}\label{prop1} Under Assumptions \textbf{A.I-II}, for every $s\geq1$ such that $BL$ does not terminate in course of the first $s$ steps, it holds $\Delta_s\leq -\rho_*$. As a result, BL algorithm terminates in at most
$$
\cM_*=32R^2/\rho_*^2+1
$$
step.
\end{proposition}
For proof, see Appendix \ref{app:prop1}.
\begin{remark}\label{rem1} Note that whatever be the origin of points $y^r$, $r\in\cR$, where $\cO$ is queried and returns separators $f_r(\cdot)$, and under Assumption \textbf{A.I}, independently of validity of Assumption \textbf{A.II}, we always have
\begin{equation}\label{eq772}
-\rho_*\geq \Delta_\cR:=\min_{y\in E_1} \max\limits_{r\in\cR}f_r(y)
\end{equation}
where $\rho_*$ is the largest of radii of Euclidean balls contained in $\cY_*\subset E_1$. In particular,  when $\Delta_{\cR}$ happens to be nonnegative (or positive), we definitely know that Assumption \textbf{A.II} does not take place (resp., that $\cY_*$ is empty).
\end{remark}
Indeed, there is nothing to prove when $\cY_*$ is empty, since then $\rho_*=-\infty$. When $\cY_*$ is nonempty, $\cY_*$ contains  a ball of radius $\rho_*\geq0$ centered at a point $y_*$; since $\|\nabla f_r(\cdot)\|_2=1$ and $f_r(\cdot)\leq0$ on $\cY_*$, we have $f_r(y_*)\leq -\rho_*$, $r\in\cR$, and (\ref{eq772}) follows. \qed
\subsection{Incorporating objective}\label{sec:bisection}

So far, our problem of interest was just a feasibility one. We can easily adjust the scheme to handle optimizing a given objective, provided that the latter is a convex function of the vector $y$ of deterministic decisions. By extending $y$, we can reduce this situation to that of minimizing a linear function $f^Ty$ with $f\neq0$ over implementable solutions $y$. This goal can be achieved by Bisection, specifically, as follows. Setting
$$
\Delta_0= \left[\min_{y\in\cY} f^Ty,\max_{y\in\cY} f^Ty\right]
$$
and selecting somehow
reliability tolerance $\delta\in(0,1)$,
optimality tolerance $\varkappa>0$, and
``stability tolerance'' $\rho>0$,
we run
 $$
 L=\Big\rfloor \log_2(|\Delta_0|/\varkappa)\Big\lfloor
  $$
  steps of Bisection, where $\rfloor a\lfloor$ is the smallest integer which is $>a$. At the $k$-th step, given the $(k-1)$-st {\sl localizer} -- segment $\Delta_{k-1}\subset\Delta_0$ of length $|\Delta_{k-1}|=2^{1-k}|\Delta_0|$, we specify our {\sl $k$-th target} $\phi_k$ as the midpoint of $\Delta_{k-1}$, add to the list of constraints specifying $\cY$ the constraint $f^Ty\leq\phi_k$ and apply to the resulting feasibility problem an algorithm of the type we have described, restricting the number of steps in this algorithm to $\cM_*(R,\rho,n)$. For example, the algorithm in question could be BL with the number of steps restricted to $M=32R^2/\rho^2+1$. Upon termination of this algorithm, the following outcomes are possible:
  \begin{enumerate}
  \item[A.] A strategic decision $y[k]\in\cY$  with $f^Ty[k]\leq\phi_k$ is found, and this decision is such that $\cO$ when queried at $y[k]$ got stuck (``productive step'')
  \item[B.] We get at our disposal a certificate of insolvability, as described in Remark \ref{rem1}, of our current feasibility problem
  \item[C.] We ran $M$ steps of the algorithm without getting stuck or running into the outcome B.
  \end{enumerate}
  When $k<L$ we pass to the next Bisection step, selecting, as our new localizer $\Delta_k$,
  \begin{itemize}
  \item the segment $\{s\in\Delta_{k-1}:s\leq \phi_k\}$ in the case of A
  \item the segment $\{s\in\Delta_{k-1}:s\geq \phi_k\}$ in the cases of B,C.
  \end{itemize}
  When $k=L$ we output, as the approximate solution $\widehat{y}$ to our optimization problem, the strategic decision found at the last productive step. Note that, by construction, the cost $f^T\bar{y}$ of decision $\widehat{y}$ does not exceed the smallest of the targets $\phi_k$ processed at productive steps $k$. If there were no productive steps, the resulting approximate  solution is undefined --- Bisection failed.
  \par
  The performance of the proposed approach can be characterized as follows. Consider the feasibility problem which we obtain when replacing the original $\cY$ with the part cut off this set by the constraint $f^Ty\leq s$. The stability number $\rho_*=\rho_*(s)$ of the resulting problem can be a positive real (the set $\cY_*[s]$ of implementable s.d.'s of the resulting feasibility problem has a nonempty interior), or zero ($\inter \cY_*[s]=\emptyset$, $\cY_*[s]\neq\emptyset$) or $-\infty$ ($\cY_*[s]=\emptyset$), and $\rho_*(s)$ is nondecreasing in $s\in\Delta_0$. Assume that
  \begin{equation}\label{eq25}
  \exists s: \rho_*(s)\geq \rho,
  \end{equation}
  and let
  $$
  s_*=\inf\{s:\rho_*(s)\geq\rho\}
  $$
  Note that $s_*\geq\underline{f}:=\min_{y\in\cY}f^Ty$, since otherwise there exists $s<\underline{f}$  with $\rho_*(s)>0$, which is impossible --  when $s<\underline{f}$, the set $\cY_*[s]$ is empty.
  \par
  Guaranteed performance of our Bisection can be described as follows:
  \begin{proposition}\label{prop2} Assume that {\rm (\ref{eq25})} takes place and, in addition,
  \begin{equation}\label{eq26}
  s_*:=\inf\{s: \rho_*(s)\geq\rho\}<\max_{y\in\cY}f^Ty-\varkappa.
  \end{equation}
  Then properly implemented Bisection\footnote{Specifically, with the sample size used by oracle to answer the $s$-th call is given by (\ref{eq20}), and the counter of oracle calls {\sl not} refreshed when passing from  a Bisection step to the next one.} with probability at least $(1-\delta)$ will output a solution $\bar{y}$ which is an $(1-\epsilon)$-implementable strategic decision with
  \begin{equation}\label{eq27}
  f^T\bar{y}\leq s_*+\varkappa
  \end{equation}
  \end{proposition}
  {\bf Proof.} It is immediately seen that with the implementation of Bisection described in the proposition the probability of the event
\begin{quote}
$\cE$: For every $k\leq L$  such that step $k$ is productive, $y[k]$ is a $(1-\epsilon)$-implementable strategic decision with $f^Ty[k]\leq\phi_k$,
and besides this, every step $k$ such that  $\phi_k>s_*$ is productive
\end{quote}
is at least $1-\delta$. Now assume that this event takes place, and let us show that in this case (\ref{eq27}) takes place. Assuming the opposite, there were no productive steps $k$ with $\phi_k\leq s_*+\varkappa$, while every step $k$ with $\phi_k>s_*$ was productive (since $\cE$ takes place). Consequently, there were no steps $k$ with $\phi_k\in\Delta:=(s_*,s_*+\varkappa]\subset\Delta_0$, so that we should either have $\Delta\subset\Delta_L$, or $\Delta\cap\Delta_L=\emptyset$. The second option can take place only if $\Delta$ is to the right of the target $\phi_k$ at some productive step $k$. We already know that the latter is not the case; thus, $\Delta\subset\Delta_L$, which is impossible, since $|\Delta_L|<\varkappa=\mes(\Delta)$, and we have arrived at a desired contradiction. \qed
\begin{remark}\label{rem2} Under Assumptions \textbf{A.I-II}, setting $f^*=\max_{y\in\cY} f^Ty$, the stability number $\rho_*>0$ as defined in
\textbf{A.II} coincides with $\rho_*(f^*)$. Denoting by $f_*$ the optimal value in the problem $\min_{y\in\cY_*} f^Ty$ of minimizing $f^Ty$ over implementable strategic decisions $y$, we clearly have $f_*<f^*$ and
$
f_*<s\leq f^*$ implying that $ \rho_*(s)\geq \rho_*{s-f_*\over f^*-f_*},
$
whence $s_*(\rho):=\inf\{s:\rho_*(s)\geq \rho\}\leq f_*+{\rho\over\rho_*}[f^*-f_*]$.
By Proposition \ref{prop2}, when using $(L+1)$-step Bisection with reliability tolerance $\delta$ and stability tolerance $\rho>0$ such that
$$
2^{-L}|\Delta_0|+{\rho\over\rho_*}[f^*-f_*]\leq f^*-f_*,
$$
the result $\bar{y}$, with probability at least $(1-\delta)$, is well defined  and is a $(1-\epsilon)$-implementable strategic decision such that
$$
f^T\bar{y}\leq f_*+2^{-L}|\Delta_0|+{\rho\over\rho_*}[f^*-f_*].
$$
\end{remark}
\section{Numerical illustration}\label{sec:numerics}

Here we present numerical results for a toy ``proof of concept'' problem, namely, an instance of the multi-product inventory problem described in Section \ref{sec:invent} where the goal is to minimize the total inventory management cost.
\paragraph{The model.} In the instance, there are $d$ products and  $K+1$ stages -- $K$ ``actual'' and one ``fictitious;'' the associated with stages blocks in strategic decision are denoted $y_1,...,y_{K+1}$. At an actual stage $t$, $1\leq t\leq K$,   $\xi_t$, $x_t$, $y_t$, and $\cZ_{\xi_t}^t$ are as described in Section \ref{sec:invent}, so that $y_t=[\ell_t;u_t;\omega_t]$, where $\ell_t\in\bR^d$ and $u_t\in\bR^d$ are the vectors of lower, resp., upper bounds on inventory level at the end of stage $t$, $\omega_t$ is an upper bound on the expenses of the stage, $x_t$ is replenishment order of the stage $t$, and $\xi_t=[d_t;o_t;h_t;p_t;r_t]$ is comprised of demand and cost/penalty coefficients observed at the stage $t$. At the fictitious stage $K+1$, $y_{k+1}\equiv \omega$ is an upper bound on the total inventory management cost, and $\xi_{K+1}\equiv\xi^K$ is the entire uncertain data.
\par
Local decisions of actual stage $t$ are replenishment orders $x_t$ of the stage, and the sets $\cZ_{\xi_t}^t$, $t\leq K$ are defined by the constraints (cf. (\ref{eq70}))
\begin{equation}\label{eq90}
\begin{array}{rl}
(a)&\ell_{t-1}+x_t-d_t\geq \ell_t,\,u_{t-1}+x_t-d_t\leq u_t,\\
(b)&o_tx_t+h_t\max[u_t,0]+p_t\max[-\ell_t,0]-r_td_t \leq \omega_t,\\
(c)&\underline{x}_t\leq x_t\leq\overline{x}_t.
\end{array}
\end{equation}
Here $\ell_0=u_0$ and $\underline{x}_t$, $\overline{x}_t$ (same as other overlined and underlined quantities below)  are part of problem's certain data. \par The local decision of the fictitious stage $K+1$ is a collection $x_{K+1}=\{\chi_1,...,\chi_K\}$ of $d$-dimensional vectors, and the set $\cZ^{K+1}_{\xi^K}$ is given by the system of constraints
\begin{equation}\label{eq92}
\begin{array}{rl}
(a)&\ell_{t-1}+\chi_t-d_t\geq \ell_t,\,u_{t-1}+\chi_t-d_t\leq u_t,\,t\leq K,\\
(b)&o_t\chi_t+h_t\max[u_t,0]+p_t\max[-\ell_t,0]-r_td_t \leq \omega_t,\,t\leq K,\\
(c)&\underline{x}_t\leq \chi_t\leq\overline{x}_t,\,t\leq K,\\
(d)&\sum_{t=1}^K[o_t\chi_t+h_t\max[u_t,0]+p_t\max[-\ell_t,0]-r_td_t]\leq \omega.
\end{array}
\end{equation}
Finally, the set $\cY$ is defined by the constraints
$$
\begin{array}{cl}
\underline{z}_t\leq \ell_t\leq u_t\leq \overline{z}_t,\,t\leq K&\hbox{[bounds on inventory levels]}\\
\underline{\omega}_t\leq \omega_t\leq \overline{\omega}_t,\,t\leq K&\hbox{[bounds on per stage budgets]}\\
\sum_{i=1}^d\max[[u_t]_i,0]\leq \overline{s},\,t\leq K&\hbox{[warehouse capacity restriction]}\\
\underline{\omega}\leq\omega\leq\overline{\omega}&\hbox{[bounds on $\omega$]}\\
\end{array}
$$
on variables $\omega$ and $u_t,\ell_t,\omega_t$, $t\leq T$, comprising a strategic decision.
\par
Note that with this formalization, implementability (or $(1-\epsilon)$-implementability) of a strategic decision  $y\in\cY$ means that for every $\xi^K\in\Xi$ (respectively, for every $\xi\in\Xi(y)$ with $P(\Xi(y))\geq1-\epsilon$)
\begin{enumerate}
\item[a)] there exist local decisions $x_t$, $t\leq K$, satisfying all constraints (\ref{eq90}) stemming from $y,\xi^K$, and thus augmenting $y$ to an implementable, uncertain data being $\xi^K$, control of the inventory, and
    \item[b)] there exists solution $x_{K+1}=\{\chi_1,...,\chi_K\}$ satisfying all constraints (\ref{eq90}) stemming from $y,\xi^K$.
     \end{enumerate}
     Assuming that $\xi^K$ and $y$ are such that a) and b) take place and
    looking at (\ref{eq90}) and (\ref{eq92}) we see that when selecting local decisions $x_t$, $t\leq K$, as  optimal solutions $\bar{x}_t$ to the (feasible !) problems of minimizing over $x_t$ the quantities $o_tx_t+h_t\max[u_t,0]+p_t\max[-\ell_t,0]-r_td_t$ under constraints (\ref{eq90})  (this selection is completely legitimate -- it specifies $\bar{x}_t$ in terms of $y$ and $\xi_t$ only) and replacing $\chi_t$ with $\bar{x}_t$,  the resulting local solution $\bar{x}_{K+1}=\{\bar{x}_t,t\leq K\}$ of the fictitious  stage $K+1$ satisfies, along with $x_{K+1}$, all stemming from $y,\xi^K$ constraints (\ref{eq92}), implying that {\sl when realization of uncertain data is $\xi^K$, the collection $(y,\bar{x}_1,\bar{x}_2,...,\bar{x}_K)$ meets all constraints of our inventory managing problem,  and  the associated total management cost
          does not exceed $\omega$}.
\par
\begin{figure}
$$
\begin{array}{ccc}
\includegraphics[scale=0.35]{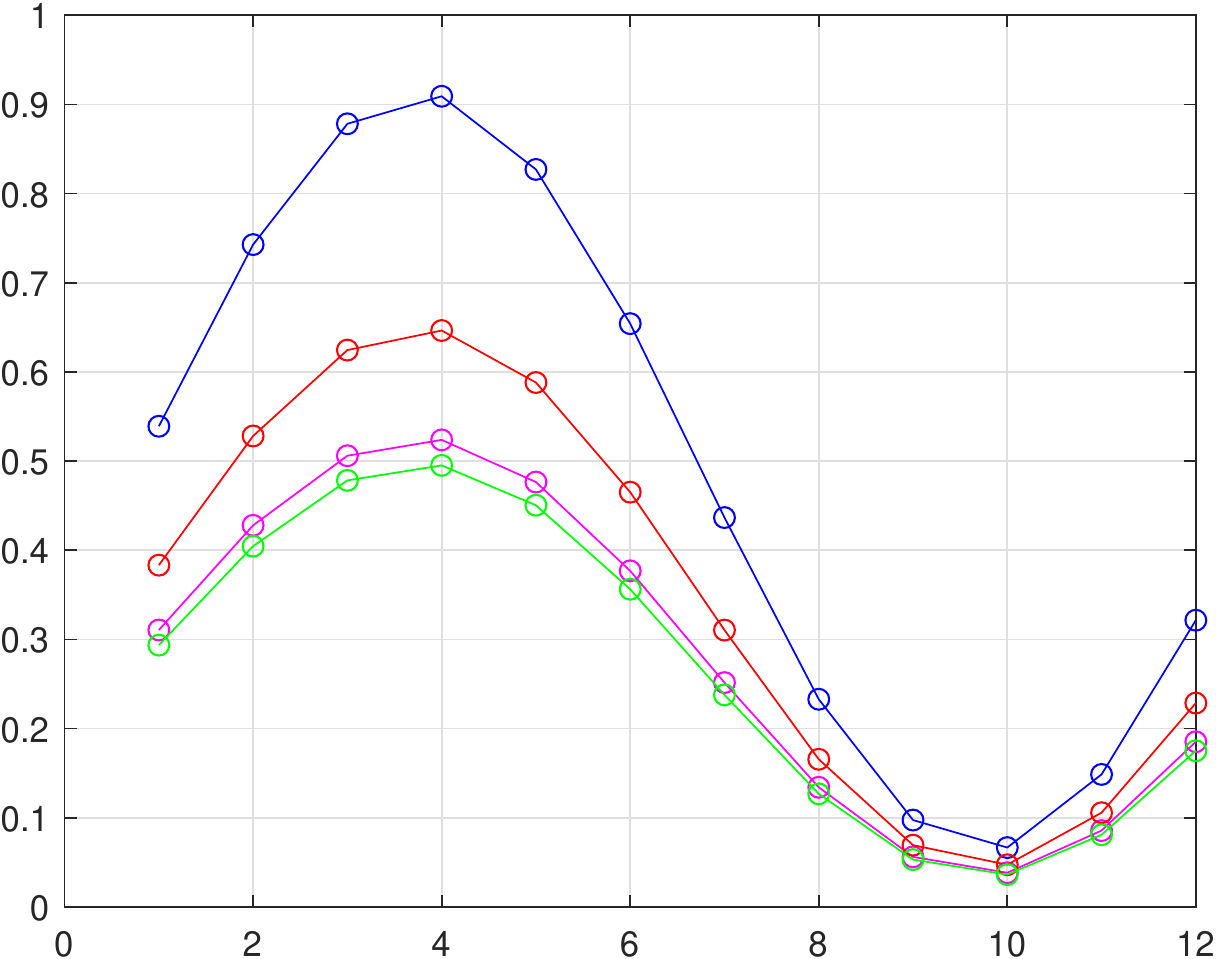}&\includegraphics[scale=0.35]{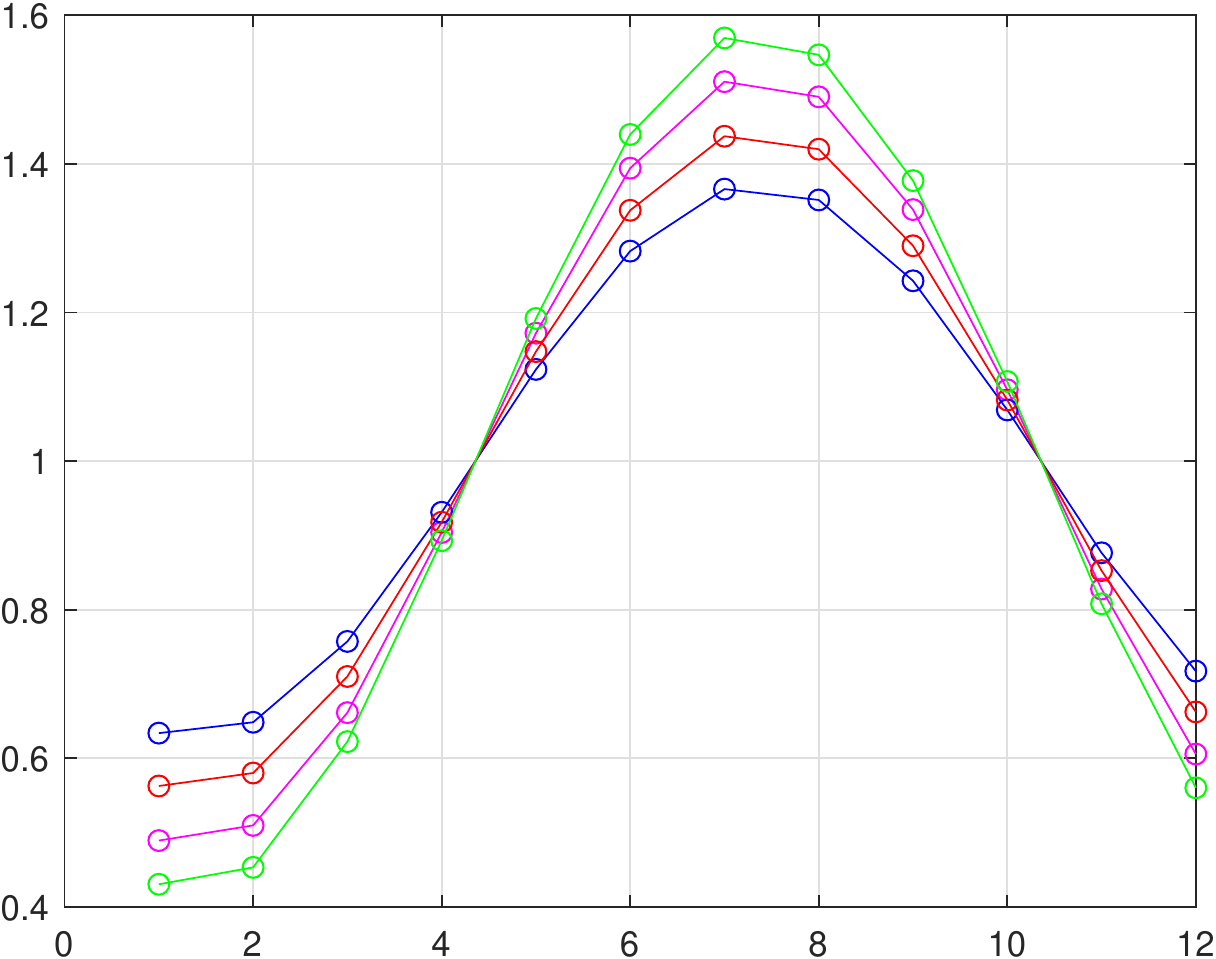}&
\includegraphics[scale=0.35]{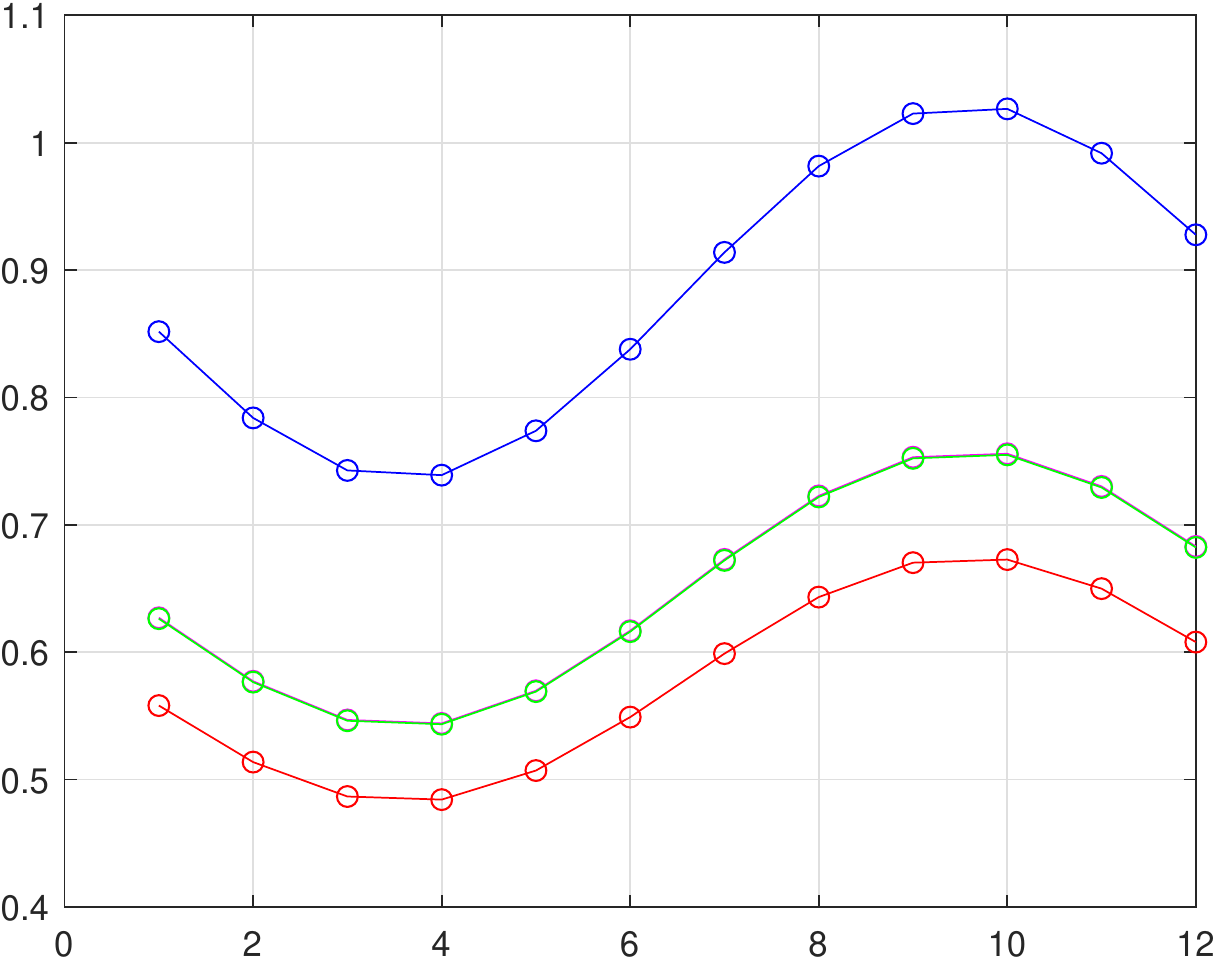}\\
\hbox{\small demands}&\hbox{\small ordering costs}&\hbox{\small holding costs}\\
\end{array}
$$
\caption{\label{Fig1} Nominal values of uncertain data vs. time}
\end{figure}

Now we can use the machinery from Section \ref{sec:bisection} to minimize the (upper bound on the) total inventory management cost $\omega$, and this is what was done in our experiments.
\paragraph{Data and results.} In the experiment we are reporting we used $d=4$, $K=12$, $\underline{z}_t=[0;0;0;0]$, $\overline{z}_t=[1;1;1;1]$.\footnote{We omit the details on how other components of certain data were specified.}
The uncertainty set $\Xi$ was the image of the set of actual uncertain data $\cH=\{\eta^K=\{\eta_t=[d_t;o_t;h_t;p_t;r_t]\}_{t\leq K}\}$ under the mapping
$\eta^K\mapsto \{\xi_t=(\eta_1,...,\eta_t)\}_{t=1}^K$, with $\cH$ and the distribution of $\eta^K$ specified as follows: each component of the uncertain $\eta$-data  --- the trajectory $d^K=(d_1,...,d_K)$ of demands and similar trajectories of ordering costs $o^K$, holding costs $h^K$, etc. ---  is uniformly distributed in its ``uncertainty box.'' For the demand, this is the box
$$
\{d^K=(d_1,...,d_K): 0.7\overline{d}_t\leq d_t \leq 1.3\overline{d}_t,\,t\leq K\}
$$
 with  positive nominal demands $\overline{d}_t$, and similarly for other components of $\eta^K$, with different components independent of each other.
Relevant\footnote{In our experiment, we used $r_t\equiv 0$; besides this, backlogged demand is forbidden due to $\underline{z}_t\equiv 0$, making the backlog penalties irrelevant.} nominal values of the uncertain data are shown in Figure \ref{Fig1}. We were looking for 0.95-implementable strategic
decision  (i.e., set $\epsilon=0.05$), and used $\delta=0.01.$\par
\begin{figure}
$$
\begin{array}{c}
\includegraphics[scale=0.9]{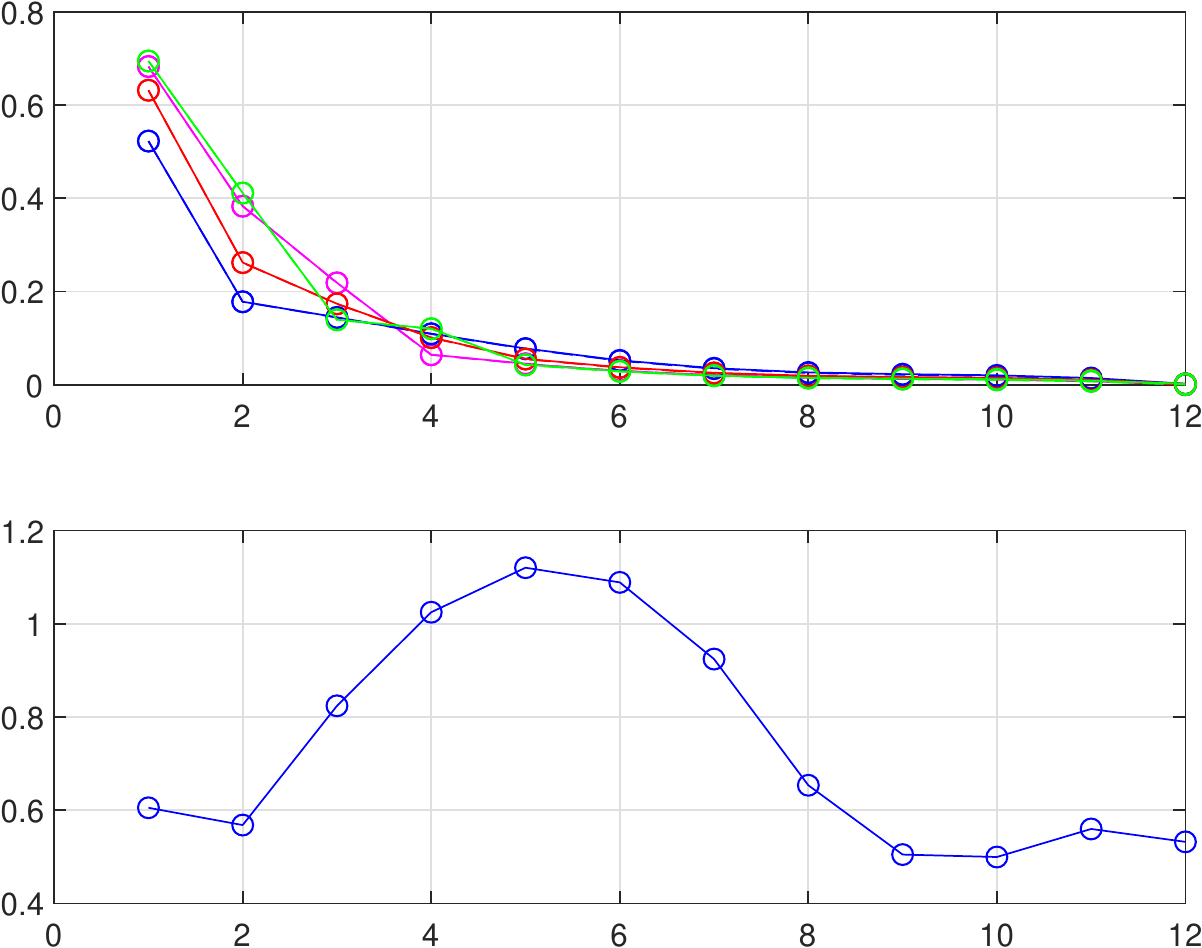}\\
\end{array}
$$
\caption{\label{Fig2} Strategic decision found by Bisection. {Top: bounds on inventory levels (the computed upper and lower bounds happened to coincide in this experiment)}; bottom: stage budgets $\omega_t$.}
\end{figure}

The numerical results obtained by processing the instance by 10-step Bisection implementing the BL algorithm are as follows. The strategic decision we got is shown in Figure \ref{Fig2}, the resulting upper bound on the inventory management cost is $\omega=19.0251$,  and the empirical management costs associated with this decision are presented in Table \ref{Tab1}. Different colors in Figures \ref{Fig1}, \ref{Fig2} correspond to different types of products operated by the inventory.
\begin{table}
$$
\begin{array}{c}
\begin{array}{||c|c|c|c||c||}
\hline
\min&\hbox{mean}&\hbox{median}&\max&\omega\\
\hline
14.7471&16.6851&16.6674&19.0245&19.0251\\
\hline
\end{array}\\
\end{array}
$$
\caption{\label{Tab1} Empirical total management cost (data over 1000 simulations) and its a priori upper bound $\omega$ as given by Bisection,
for the strategic decision shown in Figure \ref{Fig2}}
\end{table}

Note that with $0.95$-implementable strategic decision, among 1000 realizations of uncertain data $\xi^K$ we should be ready to observe about  $50$ realizations
in which we failed to augment our strategic decision by local decisions to meet all the constraints. In fact, there was just one realization of this type at all, indicating  that our construction is in fact much more reliable than is stated by our theoretical analysis.\par
To get an impression of how conservative is our decision making, we compared the associated management costs with ``utopian'' ones -- those achievable for ``clairvoyant'' decision maker who knows in advance the realization of uncertain data and selects the replenishment orders minimizing, given this realization, the total management cost.  The average, over 1000 simulations of uncertain data, excess of our management costs over the utopian ones was 10.3\%.
\section{Concluding remarks}
We want to reiterate that the only novelty, if any, in this note is  the ``tale'' about multi-stage decision making under uncertainty which allows to model risk-free decision making by the non-adjustable robust counterpart
\begin{equation}\label{eq100}
\min_y\left\{f^Ty:y\in\cY,g(y,\zeta)\leq 0\,\forall \xi\in\Xi\right\},
\end{equation} (see \cite{RO}) of uncertain convex problem with instances $\min_y \{f^Ty:y\in\cY,g(y,\zeta)\leq0\}$ parameterised by uncertain data $\zeta$ running through uncertainty set $\Xi$; here $\cY\subset\bR^n$ is a computationally tractable convex  compact set, and $g(y,\xi)$ is {\sl convex in $y$} function satisfying reasonable computability restrictions.\footnote{With the setup presented in Section \ref{sec:model}, $\zeta=\xi^K$,
$$
g(y,\xi^K)=\left\{\begin{array}{ll}+\infty,&\exists t\leq K: \{x_t:(y,x_t)\in\cZ^t_{\xi_t}\}=\emptyset,\\
0,&\hbox{otherwise},\\
\end{array}
\right.
$$
and ``reasonable computability restrictions'' reduce to computational tractability of the sets $\cZ^t_{\xi_t}$.
With slightly different setups of similar modeling power, $g(y,\zeta)$ in (\ref{eq100}) is a real-valued convex in $y$ function on $\bR^n\times\Xi$ with efficiently computable, given $y$ and $\zeta$, value and subgradient in $y$.} The remaining part is completely standard:  while semi-infinite convex problem (\ref{eq100})  can be difficult, one immediately observes that when assuming uncertain data to be stochastic: $\zeta\sim P$, and relaxing feasibility of candidate solutions to (\ref{eq100}) to $(1-\epsilon)$-feasibility  (i.e., requiring from $y\in\cY$ to satisfy $P\{\zeta:g(y,\zeta)\leq0\}\geq1-\epsilon$ instead of $g(y,\zeta)\leq0\,\forall\zeta\in\Xi=\supp(P)$), the problem, under mild additional assumptions, becomes tractable.\footnote{With the setup of Section \ref{sec:model}, ``mild additional assumptions'' require strict feasibility (nonemptiness of the interior of the feasible set) of problem (\ref{eq100}), which allows to find a near-optimal $(1-\epsilon)$-feasible solution by the machinery of  Section \ref{sec:modif}. When $g(x,\zeta)$ is real-valued  on $\bR^n\times\Xi$ and assuming (\ref{eq100}) feasible, such a solution can be found by the celebrated chance constrained approximation  originating from Calafiore and Campi \cite{CC1,CC2,CC3} and de Farias and Van Roy \cite{FVR}, see also \cite{NemShap}.}
\section{Post scriptum}
After this note was posted on arXiv, it did not take long to find out that it, basically, indeed reinvents a bicycle. The author is grateful to Prof. Jim Luedtke for pointers to two papers the author regretfully was not aware of:
\begin{itemize}
\item Vayanos, P., Kuhn, D., and Rustem, B. A constraint sampling approach for multi-stage robust optimization.  {\sl Automatica} {\bf 48:3} (2012), 459-471
\item Bodur, M. and  Luedtke, J.R. Two-stage linear decision rules for multi-stage stochastic programming. {\sl Mathematical Programming}  (2018). https://doi.org/10.1007/s10107-018-1339-4
\end{itemize}
which, for all practical purposes, cover essentially all (except, perhaps, for the BL-based algorithm for processing the model; this algorithm in any case is not the point here) the author tried to say.

\appendix
\section{Remodeling}\label{remodel}
Here we generalize the ``remodeling trick'' mentioned at the end of Section \ref{sec:invent}. Consider a general decision making model as posed in Section \ref{sec:model} and assume that
\begin{enumerate}
\item The behaviour of the controlled system on the time horizon $1,...,K$ is fully determined by realization $\xi^K$ of uncertain data and our control comprised of  decisions $y,x_1,...,x_K$. A control $y,x_1,...,x_K$ can be implemented, the uncertain data being $\xi^K$, if and only if $y\in\cY$ and $(y,x_t)\in\cZ^t_{\xi_t}$, $t\leq K$. Besides this, assume that we can split strategic decision $y$ into blocks: $y=[y_0;y_1;...;y_K]$, $y_s\in\bR^{n_s}$, in such a way that
     \begin{itemize}
     \item $y_0$ is the component of our strategic decision  which should be implemented ``at time 0,'' before the system starts to evolve, and thus the entries in $y_0$ should get numerical values when the problem is being solved, before the uncertainty reveals itself;
     \item $y_t$, $1\leq t\leq K$, is the component of our strategic decision which should be implemented at time $t$, when the components $\xi_\tau$, $\tau\leq t$, of the uncertain data are already known.
     \end{itemize}
    In this situation, we in principle could allow $y_t$ to  depend on $\xi_1,...,\xi_t$, so that our initial restriction ``strategic decision should be specified before the uncertain data starts to reveal itself'' stems from {\sl how we intend to make decisions}, and not from ``physical'' restrictions on what a decision making could be.
\item For every $t\leq K$, $\xi_t$ ``remembers'' $\xi_1,...,\xi_{t-1}$, meaning that for $s\leq t\leq K$ and $(\xi_1,...,\xi_K)\in\Xi$,\footnote{Recall that $\Xi$ is the support of the distribution $P$ of uncertain data.} $\xi_{s}$ is a deterministic function of $\xi_t$:
$$\xi_s=\Xi_{ts}(\xi_t),\,1\leq s\leq t\leq K.$$
As we remember, we can assume this w.l.o.g.
\item Sets $\cZ^t_{\xi_t}$ are of special structure:
$$
\cZ^t_{\xi_t}=\{([y_0;...;y_K],x)\in\bR^n\times\bR^{\nu_k}:([y_0;...;y_t],x)\in\cW^t_{\xi_t}\}
$$
where $\cW^t_{\xi_t}$ are given closed convex sets.
\end{enumerate}
Assume also that we are given efficiently computable basic functions $B_{sr}(\xi_s):\Xi_s\to\bR^{n_s}$, $1\leq r\leq r_s$; here $0\leq s\leq K$ and $\xi_s$ runs through $\Xi_s$, with $\Xi_0$ being a singleton.  Let us pass from deterministic strategic decisions $y$ to strategic {\sl decision rules}
\begin{equation}\label{eq300}
Y_s(\xi_s)=\sum_{r=1}^{r_s}\chi_{sr}B_{sr}(\xi_s),
\end{equation}
and treat the collection $\chi=\{\chi_{sr}:0\leq s\leq K,r\leq r_s\}$ of coefficients in (\ref{eq300}) as our new strategic decision. Given $\chi$, we control our system as follows:
\begin{itemize}
\item at time 0, we implement $y_0=Y_0$ (this rule is well defined -- $Y_0$ is a function on a singleton set $\Xi_0$ and thus is just a constant);
\item at time $t$, $1\leq t\leq K$, after $y_0,y_1,...,y_{t-1}$ and $x_0,x_1,...,x_{t-1}$ have been already built and implemented and $\xi_t$ has been observed, we compute $y_t=Y_t(\xi_t)$ according to (\ref{eq300}), find local decision $x_t$ in such a way that $([y_0;y_1;...;y_t],x_t)\in\cW^t_{\xi_t}$, and implement the decision $(y_t,x_t)$;
\item at (fictitious) time instant $K+1$, when $y=[y_0;y_1;...,y_K]$ and $x_1,...,x_K$ have been built and $\xi_{K+1}:=\xi^K$ has become known, we select fictitious local decision  $x_{K+1}=0\in\bR^{\nu_{K+1}}:=\bR$ which must satisfy the restriction $(y,x_{K+1})\in \cZ^{K+1}_{\xi^K}:=\cY\times\{0\}$.
\end{itemize}
All we need in order to {ensure that the control obtained in this way is implementable}, are the inclusions
$$
\begin{array}{ll}
(a)&\left(\chi:=\{\chi_{sr}\}_{s,r},x\right)\in\overline{\cZ}^t_{\xi_t}:=
\left\{(\chi,x):\left(\{\overline{Y}_{s}[\chi,\Xi_{ts}(\xi_t)],0\leq s\leq t\},x\right)  \in\cW^t_{\xi_t}\right\},1\leq t\leq K\\
&\hbox{where\ }\overline{Y}_{s}[\chi;\xi_s]=\sum_{r=1}^{r_s}\chi_{sr}B_{sr}(\xi_s);\\
(b)&(\chi,x_{K+1})\in \overline{\cZ}^{K+1}_{\xi^K}:=\left\{(\chi,0):
\{\overline{Y}_s[\chi;\xi_s],0\leq s\leq K\}\in\cY\right\}.\\
\\
\end{array}
$$
Indeed, inclusions $(a)$ express in terms of $\chi$-variables the restriction that in our decision making process  blocks $y_0,y_1,...,y_K$ of the strategic decision $y$ corresponding to uncertain data $\xi^K$ \footnote{These blocks are exactly vectors $\overline{Y}_s[\chi;\xi_s]=\overline{Y}_s[\chi;\Xi_{ts}(\xi_t)]$, $0\leq s\leq t\leq K$.} and local decisions $x_1,...,x_K$ meet the original dynamic constraints $(y,x_t)\in\cZ^t_{\xi_t}$, while  $(b)$  translates to the space of $\zeta$-variables the restriction $y\in\cY$. Note that  sets $\overline{\cZ}^t_{\xi_t}$ are closed and convex since $\cY$ and $\cZ^t_{\xi_t}$ are so and functions $\overline{Y}_s[\chi;\xi_s]$ are linear in $\chi$. We see that when specifying $\overline{\cY}$ as the entire space of $\chi$-variables and using this set and sets $\overline{\cZ}^t_{\xi_t}$ in the roles of $\cY$ and $\cZ^t_{\xi_t}$, our new decision making model with strategic decisions $\chi$
is of the structure described in Section \ref{sec:model} and thus can be processed in a computationally efficient fashion by the machinery we have developed. On the other hand, assuming that the set of basic functions $B_{sr}$ is reach enough to make all constant functions of $\xi_s$, $0\leq s\leq K$, linear combinations of the basic functions, our remodeling strengthens our  ``control abilities.''

\section{Proof of Proposition \ref{prop1}}\label{app:prop1}
Recall that $f_r(\cdot)$, when defined (i.e., when BL does not terminate in course of the first $r$ steps), is an affine function with $\|\nabla f_r(\cdot)\|_2=1$ which is nonpositive on $\cY_*$ and is $\geq0$ at $y^r$. Under the premise of Proposition, $\cY_*\subset E_1$ contains ball $B_*=\{y:\|y-y_*\|\leq\rho_*\}$, and since $f_r(y)\leq0$ for $y\in \cY_*$ and $\|\nabla f_r(\cdot)\|_2=1$, we have $f_r(y_*)\leq -\rho_*$. Consequently, when $r$ is such that BL does not terminate in course of the first $r$ steps, we have $f^r(y_*)\leq -\rho_*$ and thus $\Delta_r\leq-\rho_*$, as claimed.\par
Now assume that $S$ is such that BL does not terminate in course of the first $S$ steps. Observe that $f^1(\cdot)\leq f^2(\cdot)\leq...\leq f^S(\cdot)$ and therefore $\Delta_1\leq\Delta_2\leq...\leq\Delta_S\leq -\rho_*$. Setting $\delta_s=|\Delta_s|$, we get $\delta_1\geq\delta_2\geq...\geq\delta_S\geq\rho_*$. Let us split the indexes $1,...,S$ into {\sl stages} as follows. We set $s_1=S$,  $\delta^1=\delta_S$ and define the first stage $\cS_1$ as the set of all indexes $s\leq S$ such that $\delta_s\leq 2\delta_S=2\delta^1$. If $\cS_1\neq\cS:=\{1,...,S\}$, we find the largest index, $s_2$, in $\cS\backslash\cS_1$, and set $\delta^2=\delta_{s_2}$, $\cS_2=\{s\in \cS\backslash\cS_1:\delta_s\leq 2\delta^2\}$. We proceed in the same fashion: after $\delta^\ell,s_\ell,\cS_\ell$ are built, we terminate when $\cS=\cS_1\cup...\cup \cS_\ell$, otherwise select the largest index, $s_{\ell+1}$, in $\cS\backslash(\cS_1\cup...\cup\cS_\ell)$ and set $\delta^{\ell+1}=\delta_{s_{\ell+1}}$, $\cS_{\ell+1}=\{s\leq s_{\ell+1}: \delta_s\leq 2\delta^{\ell+1}\}$. Denoting by $k$ the number of the last step of this (clearly finite) process, note that
\begin{equation}\label{eq12}
\delta^{\ell+1}>2\delta^\ell,\,1\leq\ell\leq k-1 \ \&\ \delta^1\geq\rho_*
\end{equation}
Let us set
$$L_s=\{y\in E_1:f^s(y)\leq -{\half}\delta_s\},\, s\leq S,$$
 so that $y^{s+1}$ is the metric projection of $y^{s}$ onto $L_s$. We claim that for every $\ell\leq k$, all sets $L_s$, $s\in \cS_\ell$, have a point in common, specifically, the minimizer $z_\ell$ of $f^{s_\ell}$ on $E_1$. Indeed, when $s\in \cS_\ell$, we have $\delta_s\leq 2\delta^\ell$ and $f^s(\cdot)\leq f^{s_\ell}(\cdot)$, whence
$$
f^s(z_\ell)\leq f^{s_\ell}(z_\ell)=\Delta_{s_\ell}=-\delta^\ell \leq -{\half}\delta_s={\half}\Delta_s,
$$
that is, $z_\ell\in E_1$ and $f^s(z_\ell)\leq {\half}\Delta_s$, as claimed. Now, by construction $y^{s+1}\in L_s$, so that $f^s(y^{s+1})\leq -{\half}\delta_s$, $1\leq s\leq S$, and $f^s(y^s)\geq f_s(y^s)\geq0$, implying, due to evident Lipschitz continuity, with constant 1 w.r.t. $\|\cdot\|_2$,  of $f^s$, that $\|y^s-y^{s+1}\|_2\geq {\half}\delta_s$. For $s\in \cS_\ell$, $y^{s+1}$ is the metric projection of $y^s$ onto $L_s$, and $z_\ell\in L_s$, resulting in
$$\|y^{s+1}-z_\ell\|_2^2\leq \|y^s-z_\ell\|_2^2-\|y^s-y^{s+1}\|_2^2\leq \|y^s-z_\ell\|_2^2-\four\delta_s^2\leq \|y^s-z_\ell\|_2^2-\four[\delta^\ell]^2.$$
On the other hand, all points $y^s$ and sets $L_s$ belong to $E_1$, whence $\|y^{s}-z_\ell\|_2^2\leq 4R^2$ for all $s$ and $\ell$, and we conclude that
$
\Card(\cS_\ell)\leq 16R^2/[\delta^\ell]^2,\,1\leq\ell\leq k.
$
Invoking (\ref{eq12}), we arrive at
$$
S=\sum_{\ell=1}^k\Card(\cS_\ell)\leq 16R^2\sum_{\ell=1}^k [\delta^\ell]^{-2}\leq 16R^2\sum_{\ell=1}^k [\delta^1]^{-2}2^{-4(\ell-1)}
\leq 32R^2/\rho_*^2.\eqno{\hbox{\qed}}
$$

\end{document}